\documentclass[nummat,final,numbook]{svjour}
\usepackage{psfig}
\input epsf

\def\pref#1{(\ref{#1})}
\def\dsp{\displaystyle}
\def\Frac#1#2{\frac
{
 {\raise.6ex
 \hbox{$\displaystyle#1$}}
}
{
 {\lower.6ex
 \hbox{$\displaystyle#2$}}
 }
}

\def\begeq{\begin{equation}
\begin{array}{ll}
}
\def\endeq{\end{array}
\end{equation}}

\def\aof{{a^{\frac14}}}
\def\atf{{a^{\frac34}}}
\def\intp{\int_0^\infty}

\def\bo{{\cal O}}
\def\C{{{\cal C}}}
\def\ph{{\rm ph}}
\def\Frac#1#2{\frac{\displaystyle{#1}}{\displaystyle{#2}}}
\def\wt{\widetilde}
\def\tfrac#1#2{{{\lower.6ex
\hbox{$\scriptstyle#1$}}\over 
{\raise.7ex
\hbox{$\scriptstyle#2$}}}}

\def\RR{{{\rm I}\!{\rm R}}}

\begin{document}
\title{Integral Representations for Computing \\
Real Parabolic Cylinder Functions}
\author{Amparo Gil\inst{1} \and Javier Segura\inst{2}
 \and Nico M. Temme\inst{3}
}                     
%
%
\institute{
Departamento de Matem\'aticas, U. Aut\'onoma de Madrid, 28049-Madrid, Spain; 
amparo.gil@uam.es
 \and Departamento de Matem\'aticas, Estad\'{\i}stica y 
        Computaci\'on,
        U. Cantabria, 39005-Santander, Spain; javier.segura@unican.es
\and CWI, P.O. Box 94079, 1090 GB Amsterdam, The Netherlands; nicot@cwi.nl }
\date{Received: date / Revised version: date}
%
\maketitle
\begin{abstract}

Integral representations are derived for the parabolic cylinder functions
$U(a,x)$, $V(a,x)$ and $W(a,x)$ and their derivatives.  The new integrals
will be used in numerical algorithms based on quadrature.  They follow from
contour integrals in the complex plane, by using methods from asymptotic
analysis (saddle point and steepest descent methods), and are stable
starting points for evaluating the functions $U(a,x)$, $V(a,x)$ and
$W(a,x)$ and their derivatives by quadrature rules.  In particular, the new
representations can be used for large parameter cases.  Relations of the
integral representations with uniform asymptotic expansions are also given.
The algorithms will be given in a future paper.
\end{abstract}

\newpage
\section*{Contents of the paper}\label{sec:I.OV}
We give an overview of the structure of the paper.
\begin{description}
  \item[Section 1\ ] The basic properties of the parabolic 
                     cylinder functions $U(a,z)$ and $V(a,z)$ 
                     that are used in this paper. 
                      
  \item[Section 2\ ] The integral representation of $U(a,x)$ for $a>0$.
    \begin{description}
    \item[2.1\ ] the case $x\ge0$. 
    \item[2.2\ ] the case $x\le0$. 
    \item[2.3\ ] a Wronskian relation for 4 integrals.
    \item[2.4\ ] the relation with uniform asymptotic expansions.
      \end{description}
 
 \item[Section 3\ ] The integrals of $U(a,x)$ and $V(a,x)$ for $a<0$.
    \begin{description}
    \item[3.1\ ] the case $-1\le t \le 1$, where $t=x/(2\sqrt{|a|})$. 
      \begin{description}
        \item[3.1.1\ ] a Wronskian relation for 4 integrals. 
        \item[3.1.2\ ] the relation with uniform asymptotic expansions. 
      \end{description}

    \item[3.2\ ] the case $t\ge1$.
      \begin{description}
        \item[3.2.1\ ] a Wronskian relation for 4 integrals.
        \item[3.2.2\ ] the relation with uniform asymptotic expansions. 
      \end{description}
    \item[3.3\ ] the case $t\sim1$.
    \item[3.4\ ] the case $t\le -1$.
    \end{description}
 
\item[Section 4\ ]  The $W-$function.    
  \begin{description}
    \item[4.1\ ] the standard solutions.
    \begin{description}
        \item[4.1.1\ ] the function $\rho(a)$.
      \end{description}
    \item[4.2\ ] the case $a<0$.
    \item[4.3\ ] the case $a>0$.
    \begin{description}
        \item[4.3.1\ ] the case $t\ge1$, where $t=x/(2\sqrt{a})$.
        \item[4.3.2\ ] the case $-1\le t\le 1$. 
        \item[4.3.3\ ] unstable representations. 
      \end{description}  
     \end{description}

\item[Section 5\ ] Concluding remarks.    
 \end{description}
\section{Introduction}
\label{sec:int}
The solutions of the differential equation
\begin{equation} \label{int:i1}
\frac{d^2y}{dz^2}-\left(\tfrac14z^2+a\right)y=0
\end{equation}
are called parabolic cylinder functions and are entire functions of $z$. 
As in \cite{abst}, Chapter 19, \cite{olver}, and \cite{temsf} we denote two
standard solutions of \pref{int:i1} by $U(a,z), V(a,z)$.  Another notation
found in the literature is $D_\nu(z)=U(-\nu-\tfrac12,z)$.  Special cases
are Hermite polynomials, error functions and Fresnel integrals.

Values at the origin are given by
\begin{equation}\label{int:i2}
\begin{array}{lll}
\quad
U(a,0)=\frac{\sqrt{\pi}}{2^{\frac12 a+\frac14}\Gamma(\frac34+\frac12a)}&, &
U'(a,0)=-\frac{\sqrt{\pi}}{2^{\frac12a-\frac14}\Gamma(\frac14+\frac12a)},\\
\\
V(a,0)=\frac{\pi\ 2^{\frac12a+\frac14}}
{[\Gamma(\frac34-\frac12a)]^2\Gamma(\frac14+\frac12a)}&, &
V'(a,0)=\frac{\pi\ 2^{\frac12a+\frac34}}
{[\Gamma(\frac14-\frac12a)]^2\Gamma(\frac34+\frac12a)}.
\end{array}
\end{equation}
Then we have
\begin{equation}\label{int:i4}
U(a,z)=U(a,0)\,y_{1}(a,z)+U'(a,0)\,y_{2}(a,z),  
\end{equation}
\begin{equation}\label{int:i5}
V(a,z)=V(a,0)\,y_{1}(a,z)+V'(a,0)\,y_{2}(a,z),  
\end{equation}
where
\begin{equation}
\begin{array}{ll} \label{int:i6}
y_1(a,z)
& = e^{\tfrac14z^2}{}_1F_1\left(-\tfrac12a+\tfrac14,\tfrac12;-\tfrac12z^2\right)
\\ 
& =e^{-\tfrac14z^2}{}_1F_1\left(\tfrac12a+\tfrac14,\tfrac12;\tfrac12z^2\right),
 \\
y_2(a,z)
& = ze^{\tfrac14z^2}{}_1F_1\left(-\tfrac12a+\tfrac34,\tfrac32;-\tfrac12z^2\right)\\
& =ze^{-\tfrac14z^2}{}_1F_1\left(\tfrac12a+\tfrac34,\tfrac32;\tfrac12z^2\right),
\end{array}
\end{equation}
and the confluent hypergeometric function is defined by
\begin{equation} \label{int:i7}
{}_1F_1(a,c;z)=\sum_{n=0}^\infty \frac{(a)_n}{(c)_n}\,\frac{z^n}{n!},
\end{equation}
with $(a)_n=\Gamma(a+n)/\Gamma(a), n=0,1,2,\ldots$.

The functions $y_1(a,z)$ and $y_2(a,z)$ are the simplest even and odd solutions
of \pref{int:i1} and the Wronskian of this pair is given by
\begin{equation}\label{int:i8}
\begin{array}{ll} 
{{\cal W}}[y_1(z),y_2(z)]=y_1(z)y_2'(z)-y_1'(z)y_2(z)=1.
\end{array}
\end{equation}
>From a numerical point of
view, the pair $\{y_1,y_2\}$ is not a satisfactory pair \cite{mil52}, 
because they have almost the same asymptotic behaviour at infinity.

The behaviour of $U(a,z)$ and  $V(a,z)$ is, for large positive $z$ and $z\gg|a|$:
\begin{equation}\label{int:i9}
\begin{array}{ll} 
U(a,z)&=e^{-\tfrac14z^2}z^{-a-\tfrac12}\left[1+\bo\left(z^{-2}\right)\right],\\
V(a,z)&=\sqrt{{2/\pi}}e^{\tfrac14z^2}z^{a-\tfrac12}\left[1+\bo\left(z^{-2}\right
)\right].
\end{array}
\end{equation}
Clearly, numerical computations of $U(a,z)$ that are based on the 
representations
in \pref{int:i4} and \pref{int:i5} should be done with great care, because of 
the 
loss of accuracy if $z$ becomes large. Also, for large $a$ these 
representations become useless.

The Wronskian relation between $U(a,z)$ and $V(a,z)$ reads:
\begin{equation} \label{int:i10}
{{\cal W}}[U(a,z),V(a,z)]  =
 \sqrt{2/\pi}.   
\end{equation}
\begin{equation} \label{int:i11}
{{\cal W}}[U(a,z),U(a,-z)]=
 \frac{\sqrt{2\pi}}{\Gamma(a+\frac12)}.   
\end{equation}
which shows that $U(a,z)$ and $V(a,z)$ are independent 
solutions of \pref{int:i1}
for all values of  $a$. Other relations are
\begin{equation}\label{int:i12}
\begin{array}{ll} 
U(a,z)&=\Frac{\pi}{\cos^2\pi a\,\Gamma(a+\tfrac12)}
\left[V(a,-z)-\sin\pi a\,V(a,z)\right],\\
V(a,z)&=\Frac{\Gamma(a+\tfrac12)}\pi
\left[\sin\pi a\,U(a,z)+U(a,-z)\right].
\end{array}
\end{equation}

Equation \pref{int:i1} has two turning points at $\pm 2\sqrt{{-a}}$.  For
real parameters they become important if $a$ is negative, and the
asymptotic behaviour of the solutions of \pref{int:i1} as $a\to-\infty$
changes significantly if $z$ crosses the turning points.  At these points
Airy functions are needed for describing the asymptotic behaviour.

The purpose of this paper is to give integral representations of $U(a,x)$
and $V(a,x)$ for real values of $a$ and $x$.  We use integral
representations from the literature and modify these by saddle point
methods.  In this way we obtain integrands that are non-oscillating, also
for the case $a<0$.  In particular, we can use the new representations for
large parameter cases.   In earlier papers \cite{temsteep} and \cite{gilscorer}
we have used these methods for obtaining stable integral representations for 
modified Bessel functions with pure imaginary order and for 
inhomogeneous Airy functions (Scorer functions).

We give relations of the
integral representations with uniform asymptotic expansions,
which are taken from \cite{olpar} and \cite{tempar}. 
We only give the expansions in terms of elementary functions. Uniform 
expansions in terms of Airy functions can be found in \cite{olpar}, 
and a modified form in \cite{tempar}.

We also consider solutions $W(a,\pm x)$ of 
the differential equation 
\begin{equation} \label{int:i13}
W''+\left(\tfrac14x^2-a\right)W=0,
\end{equation}
a modified form of 
\pref{int:i1}, again for real $a$ and $x$. Properties of $W(a,x)$ are given 
in \S4, which can be found in \cite{abst} and \cite{mil55}.

In a future paper we give algorithms based on quadrature rules for
evaluating the integral representations of $U(a,x)$, $V(a,x)$ and $W(a,x)$.

In \cite{tempar} numerical and asymptotic aspects of the parabolic
cylinder functions have been discussed, and we refer to this paper 
frequently.  
The notation of certain quantities is also as in  \cite{tempar}.  
The asymptotic methods referred to
in this paper (saddle point methods) can be found in \cite{olver} and 
\cite{wong}. For an overview of the numerical aspects and software for the 
parabolic cylinder 
functions we refer to \cite{lozol}.

\section{Integral representations  for {\protect\boldmath $a>0$}}
\label{sec:uvp}%
We derive integral representations for $U(a,x)$ and $U(a,-x)$. The computation
of $V(a,x)$ for $a>0$ can be based on the second relation in \pref{int:i12}. 
For $a>0$ the functions $U(a,x)$ and $U(a,-x)$ have a non-vanishing
Wronskian relation (see \pref{int:i11}), and moreover, these functions 
constitute
a numerically satisfactory pair of solutions of \pref{int:i1}.

\subsection{The case {\protect\boldmath $x\ge0$}}
\label{subsec:uvpp}%
We take the integral (see\cite{abst}, formula 19.5.4)
\begin{equation}\label{pp:int}
U(a,x)=\frac{e^{\frac14x^2}}{i\sqrt{2\pi}}\int_{\C} e^{-xs+\frac12s^2}s^{-a}
\frac{ds}{\sqrt{s}},
\end{equation}
where $\C$ is a vertical line on which $\Re{s}>0$. 
On $\C$ we have $-\frac12\pi<\ph{\,s}<\frac12\pi$,
and the many-valued function $s^{-a-1/2}$ assumes its principal value.
The transformations
\begin{equation}\label{pp:tr}
x=2t\sqrt{a},\quad s=\sqrt{a}\,w
\end{equation}
give
\begin{equation}\label{pp:intw}
U(a,x)=\frac{e^{\frac14x^2}a^{\frac14-\frac12 a}}{i\sqrt{2\pi}}\int_{\C} 
e^{a\phi(w)}
\frac{dw}{\sqrt{w}},
\end{equation}
where
\begin{equation}\label{pp:phi}
\phi(w)=\tfrac12w^2-2tw-\ln w.
\end{equation}
The saddle points follow from solving
\begin{equation}\label{pp:phip}
\phi'(w)=\frac{w^2-2tw-1}{w}=0,
\end{equation}
giving saddle points at $t\pm\sqrt{t^2+1}$. We take for the path $\C$ in
\pref{pp:intw} the vertical through the positive saddle point
\begin{equation}\label{pp:w0} 
w_0=t+\sqrt{t^2+1}.
\end{equation}
At this saddle point $\C$ coincides with the steepest descent path
trough $w_0$.  The complete steepest descent path follows from solving 
$\Im[\phi(w)]=\Im[\phi(w_0)]$. 
In the present case $\Im[\phi(w_0)]=0$ and we obtain for the saddle 
point contour the equation
\begin{equation}\label{pp:spc}
\tfrac12r^2\sin2\theta-2tr\sin\theta-\theta=0, 
\quad{\rm where}\quad w=re^{i\theta},
\end{equation}
which can be solved for $r=r(\theta)$:
\begin{equation}\label{pp:r}
r=\frac{t+\sqrt{t^2+\theta\cot\theta}}{\cos\theta},
\quad -\tfrac12\pi<\theta<\tfrac12\pi.
\end{equation}

\vspace*{0.3cm}

\begin{center}
\epsfxsize=4cm \epsfbox{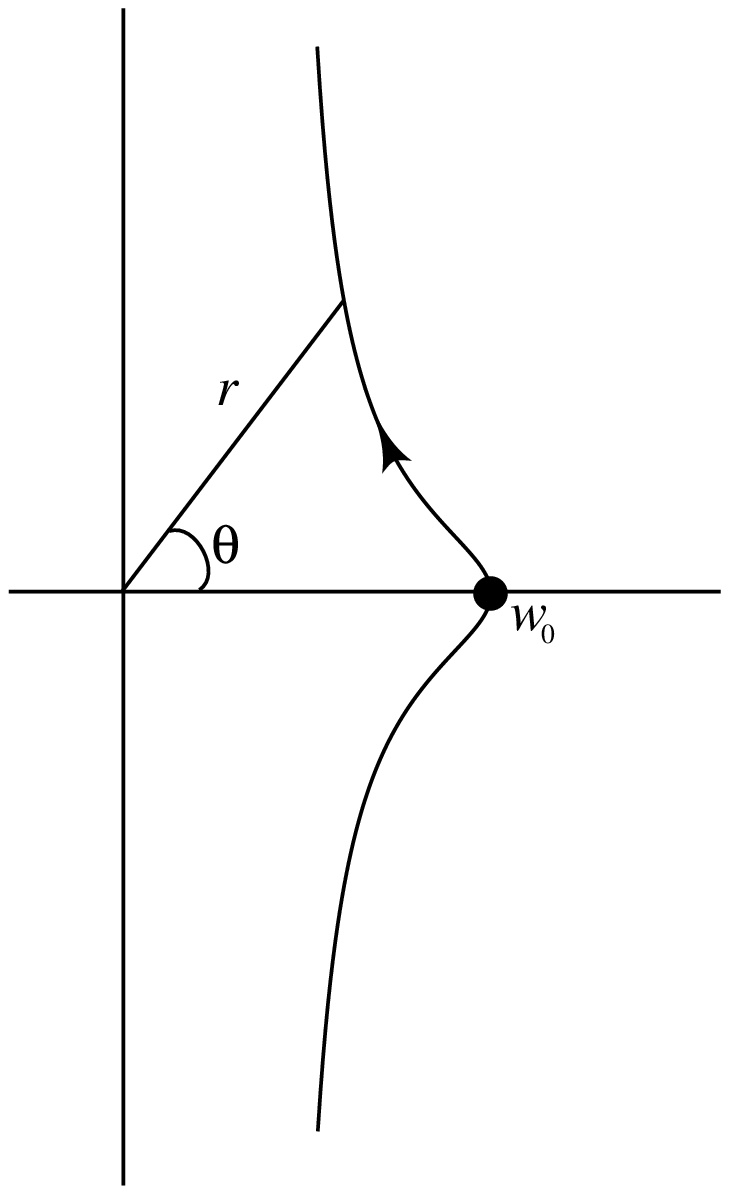}
\end{center}

\begin{quote}
{\bf Figure 1.}\quad
Steepest descent contour for the integral in \pref{pp:intw}.
\end{quote}

Then \pref{pp:intw} can be written as
\begin{equation}\label{pp:inttheta}
U(a,x)=\frac{e^{\frac14x^2+a\phi(w_0)}a^{\frac14-\frac12 a}}{\sqrt{2\pi}}
\int_{-\frac12\pi} ^{\frac12\pi}
e^{a\psi(\theta)} g(\theta)\,d\theta,
\end{equation}
where
\begin{equation}\label{pp:psi}
\psi(\theta)=\Re[\phi(w)-\phi(w_0)]=\tfrac12r^2\cos2\theta-2tr\cos\theta
-\ln r-\phi(w_0),
\end{equation}
and
\begin{equation}\label{pp:gt}
\begin{array}{ll}
g(\theta)& =\Im\left[\frac1{\sqrt{w}}\frac{dw}{d\theta}\right]
=\Im\left[e^{\frac12i\theta}\frac1{\sqrt{r}}
\left(\frac{dr}{d\theta}+ir\right)\right]\\
\\
& =
\Frac{(2\cos\theta+1)r^2-2tr+1}
{4\sqrt{r}\cos\frac12\theta\sqrt{t^2+\theta\cot\theta}}.
\end{array}
\end{equation}
The quantity $\wt\xi$ defined by
\begin{equation}\label{pp:xit}
\wt\xi=\tfrac12[t\sqrt{t^2+1}+\ln(t+\sqrt{t^2+1})],
\end{equation}
is used in the asymptotic representation of $U(a,x)$ in this case; 
see \cite{tempar}, formula (2.29). We have
\begin{equation}\label{pp:xitt}
\tfrac14x^2+a\phi(w_0)=a[\tfrac12-t\sqrt{t^2+1}-
\ln(t+\sqrt{t^2+1})]=a(\tfrac12-2\wt\xi).
\end{equation}

This gives
\begin{equation}\label{pp:intas}
U(a,x)=\frac{\aof\,e^{-2a\wt\xi}}{\sqrt{2\pi} \gamma(a)}
\int_{-\frac12\pi} ^{\frac12\pi}
e^{a\psi(\theta)} g(\theta)\,d\theta,
\end{equation}
where 
\begin{equation}\label{pp:gammaa}
\gamma(a)= e^{-\frac12a} a^{\frac12a}.
\end{equation}

For the derivative $U'(a,x)$ we can start from
\pref{pp:int}, and we have
\begin{equation}\label{pp:intd}
U'(a,x)=\frac{e^{\frac14x^2}}{i\sqrt{2\pi}}
\int_{\C} e^{-xs+\frac12s^2}s^{-a}
(\tfrac12x-s)\frac{ds}{\sqrt{s}},
\end{equation}
This can be written as 
\begin{equation}\label{pp:intthetad}
U'(a,x)=-\frac{\atf\,e^{-2a\wt\xi}}{\sqrt{2\pi} \gamma(a)}
\int_{-\frac12\pi} ^{\frac12\pi}
e^{a\psi(\theta)} h(\theta)\,d\theta,
\end{equation}
where
\begin{equation}\label{pp:ht}
\begin{array}{ll}
h(\theta)&{\dsp
=-\Im\left[\frac1{\sqrt{w}}\frac{dw}{d\theta}(t-w)\right]
}\\
&\\
&{\dsp
=\frac{r^3-tr^2(2\cos\theta-1)+r(2t^2+1+2\cos\theta)-t}
{4\sqrt{r}\cos\frac12\theta\sqrt{t^2+\theta\cot\theta}}.
}
\end{array}
\end{equation}

\subsection{The case {\protect\boldmath $x\le0$}}
\label{subsec:uvpn}%
This case can be done by using the representation of the previous section.
However, when $t$ is a large negative number, the saddle point $w_0$ defined in
\pref{pp:w0} is close to origin, at which point the integrand of \pref{pp:intw}
is singular. As a consequence, the functions $\psi(\theta)$ and $g(\theta)$ in
\pref{pp:inttheta} have singularities close to the origin $\theta=0$ when $t$ is
a large negative number.

In the present case we start with the well-known integral representation
(see \cite{abst}, formula 19.5.3)
\begin{equation}\label{pn:ints}
U(a,-x)=\frac{e^{-\tfrac14x^2}}{\Gamma(a+\tfrac12)}
\intp s^{a-\tfrac12} e^{-\tfrac12s^2+xs}\,ds,\quad a>-\tfrac12.
\end{equation}
There are no oscillations, but it is convenient to transform the 
integral in such a way that 
the saddle point is at the origin and a suitable normalization is obtained.
The transformations \pref{pp:tr} give
\begin{equation}\label{pn:intw}
U(a,-x)=\frac{a^{\frac12a+\frac14}\,e^{-\tfrac14x^2}}{\Gamma(a+\tfrac12)}\intp
e^{-a\phi(w)}\, \frac{dw}{\sqrt{w}},
\end{equation}
where $\phi(w)$ is given in \pref{pp:phi}. 
The positive saddle point $w_0$ is as in \pref{pp:w0}. We transform this point
to the origin by writing $w=w_0(1+u)$, which gives
\begin{equation}\label{pn:intrho} 
U(a,-x)=\frac{\aof\,\sqrt{w_0}\,\gamma(a)\,e^{2a\wt\xi}}
{\Gamma(a+\tfrac12)}\int_{-1}^{\infty} e^{-a\psi(u)}\,
\frac{du}{\sqrt{1+u}},
\end{equation}
where we have used \pref{pp:xitt}, $\gamma(a)$ is defined in \pref{pp:gammaa}, 
and 
\begin{equation}\label{pn:psiu}
\psi(u)=\tfrac12 w_0^2 u^2+u-\ln(1+u).
\end{equation}
For the derivative we have
\begin{equation}\label{pn:intrhod} 
\begin{array}{ll}
U'(a,-x)=&\dsp{-\frac{\atf\,\sqrt{ w_0}\, \gamma(a) \,e^{2a\wt\xi}}
{\Gamma(a+\tfrac12)}}\\
&\dsp{\times \int_{-1}^{\infty} e^{-a\psi(u)}\,
\left(\sqrt{t^2+1}+w_0u\right)
\frac{du}{\sqrt{1+u}}},
\end{array}
\end{equation}

To avoid numerical cancellation for small values of $u$
in the computation of $\psi(u)$ defined in \pref{pn:psiu},
a specific code is needed for the evaluation of $\ln(1+u)-u$. 

\subsection{A Wronskian for the integrals}
\label{subsec:uvppwr}%
When checking the numerical algorithms the Wronskian relations in 
\pref{int:i10} and \pref{int:i11} can be used. When the parameters are 
large it is more convenient to use a Wronskian relation that is based on the 
integrals derived in the section. This gives a better control of the errors 
that occur in the quadrature rules, because large and small factors are 
not present in the integrals.

We write (see \pref{pp:intas}, \pref{pp:intthetad}, 
\pref{pn:intrho},
and \pref{pn:intrhod}, respectively)

\begin{equation}\label{pp:inti}
U(a,x)=\frac{\aof\,e^{-2a\wt\xi}}{\sqrt{2\pi}\, \gamma(a)}\,
I(a,x),
\end{equation}

\begin{equation} \label{pp:intj} 
U'(a,x)=-\frac{\atf\,e^{-2a\wt\xi}}{\sqrt{2\pi}\,\gamma(a)}\,
I_d(a,x),
\end{equation}

\begin{equation} \label{pn:intid} 
U(a,-x)=
\frac{\aof\sqrt{w_0}\gamma(a)\,e^{2a\wt\xi}}
{\Gamma(a+\tfrac12)}\,
J(a,x),
\end{equation}

\begin{equation} \label{pn:intjd} 
U'(a,-x)=
-\frac{\atf\sqrt{w_0}\gamma(a)\,e^{2a\wt\xi}}
{\Gamma(a+\tfrac12)}\,
J_d(a,x).
\end{equation}
Then the relation for the integrals reads
\begin{equation} \label{ppn:wr} 
I(a,x)\,J_d(a,x)+I_d(a,x)\,J(a,x) = \frac{2\pi}{a \sqrt{w_0}}.
\end{equation}

\subsection{Uniform asymptotic expansions  for {\protect\boldmath $a>0$}}
\label{subsec:uvppas}%
The quantities related with the integrals are closely related with the 
uniform asymptotic expansions given in formulas (2.29), 
(2.33) and (2.34) of \cite{tempar}. We have
\begin{equation} \label{ppn:as1}
\begin{array}{ll} 
I(a,x)=\frac{\sqrt{\pi}}{\sqrt{a} (t^2+1)^{1/4}}\,\wt{F}_\mu(t),&
I_d(a,x)=\frac{\sqrt{\pi}(t^2+1)^{1/4}}{\sqrt{a} }\,\wt{G}_\mu(t),
\\
J(a,x)=\frac{\sqrt{\pi }}{\sqrt{aw_0} (t^2+1)^{1/4}}\,\wt{P}_\mu(t),&
J_d(a,x)=\frac{\sqrt{\pi }(t^2+1)^{1/4}}{\sqrt{aw_0} }\,\wt{Q}_\mu(t),
\end{array}
\end{equation}
where $\wt{F}_\mu(t)$, $\wt{G}_\mu(t)$, $\wt{P}_\mu(t)$, and $\wt{Q}_\mu(t)$ are 
supplied with asymptotic expansions that have a double asymptotic property: 
one of the parameters $a$ or $t$ (or both) should be large. Recurrence 
relations for the coefficients of the expansions are given in \cite{tempar}.

\section{Integral representations  for {\protect\boldmath $a<0$}}
\label{sec:uvn}%
We give integral representations for $U(-a,x)$ and $V(-a,x)$, with $a>0$,
and we consider three $x-$intervals. Let $t=x/(2\sqrt{a})$. 
The differential equation 
\pref{int:i1} becomes for $U(-a,2t\sqrt{a})$ and $V(-a,2t\sqrt{a})$
in terms of $t$
\begin{equation} \label{n:de}
\frac{d^2y}{dt^2}-4a^2\left(t^2-1\right)y=0,
\end{equation}
which has turning points at $t=\pm 1$. Consequently, we
consider the intervals $t\le-1$, $|t|\le 1$ and $t\ge 1$. 
We start with the middle interval, where the oscillations occur.

\subsection{The case {\protect\boldmath $-1\le t\le 1$}}
\label{subsec:uvnm}%
We consider the integral 
\begin{equation}\label{nm:ints}
Y(a,x)=
\int_0^{\infty} e^{-\frac12s^2+xis}s^{a-\frac12}\,ds,\quad \Re a > -\tfrac12.
 \end{equation}
Using \pref{pn:ints}, we see that
\begin{equation}\label{nm:yu}
Y(a,x)= \Gamma (a+\tfrac12 )
e^{-\frac14x^2}\, U(a,-ix).
\end{equation}
We also have
\begin{equation}\label{nm:uuv}
\sqrt{\tfrac2\pi}e^{-\frac12\pi ia+\frac14\pi i}
U(a,-ix)=
U(-a,x)/\Gamma (a+\tfrac12 )+i V(-a,x).
\end{equation}
This follows from using the initial values in \pref{int:i2}
and those of $Y(a,x)$. It also follows from the relations in \pref{int:i12}
and 19.4.6 in \cite{abst}. 

Hence, 
\begin{equation}\label{nm:yuv}
\sqrt{\tfrac2\pi}e^{-\frac12\pi ia+\frac14\pi i} e^{\frac14x^2} Y(a,x)=
U(-a,x)+i\Gamma (a+\tfrac12 )V(-a,x).
 \end{equation}
We see that the single integral \pref{nm:ints}
produces $U(-a,x)$ and $V(-a,x)$ by taking real and imaginary parts.

We proceed with $Y(a,x)$, and the transformations
$x=2\sqrt{a}t$, $s=\sqrt{a}w$ give
\begin{equation}\label{nm:intw}
Y(a,x)= a^{a/2+1/4}
\int_0^{\infty} e^{-a\phi(w)}\frac{dw}{\sqrt{w}},
 \end{equation}
where 
\begin{equation}\label{nm:phiw}
\phi(w)=\tfrac12w^2-2itw-\ln w.
 \end{equation}
We consider a path through the saddle point
\begin{equation}\label{nm:wp}
w_+=it+\sqrt{1-t^2}.
\end{equation}
We have
\begin{equation}\label{nm:phiwp}
\phi(w_+)=\tfrac12+t^2+2i\left(\eta-\tfrac14\pi\right),\quad 
\eta = \tfrac12\left(\arccos t-t\sqrt{1-t^2}\right),
 \end{equation}
where $\arccos t$ has values in $[0,\pi]$ for $t\in[-1,1]$.

The path of steepest descent starts at  $w=0$, runs through $w_+$,
and terminates at $+\infty$; see Figure 2. 
The path follows from solving the equation
\begin{equation}\label{nm:imphi}
\Im\phi(w)=\Im\phi(w_{+}),
\end{equation}
that is, from solving
\begin{equation}\label{nm:rt}
\tfrac12r^2\sin2\theta-2tr\cos\theta-\theta-2\eta+\tfrac12\pi=0,
 \end{equation}
where $w=re^{i\theta}$.
The solution of \pref{nm:rt} reads
\begin{equation}\label{nm:r}
\begin{array}{l}
r=\dsp{\frac{t\cos\theta+
\sigma\sqrt{t^2\cos^2\theta+\sin\theta\cos\theta(\theta+2\eta-\frac12\pi)}}
{\sin\theta\cos\theta}}, \\
\\
0\le\theta\le\theta_0,
\end{array}
\end{equation}
where $\theta_0=-2\eta+\frac12\pi$; the square root is non-negative. The 
number $\sigma$ equals $-1$ when $\ph\, w_+\le\theta\le\theta_0$, and 
$+1$ when 
$0\le \theta\le\ph\, w_+$. Observe that $\ph\, w_+=\frac12\pi-\arccos t$. When
$\theta=\theta_0$, we have $r=0$; 
when $\theta=\ph\, w_+$, we have $r=1$,
and $\theta=0$ gives $r=\infty$. 
For $t=0$ the path coincides with the positive real axis.
When $t<0$ the paths are in the lower half plane, and follow 
from those for $t>0$ by symmetry.

\begin{center}
\epsfxsize=10cm \epsfbox{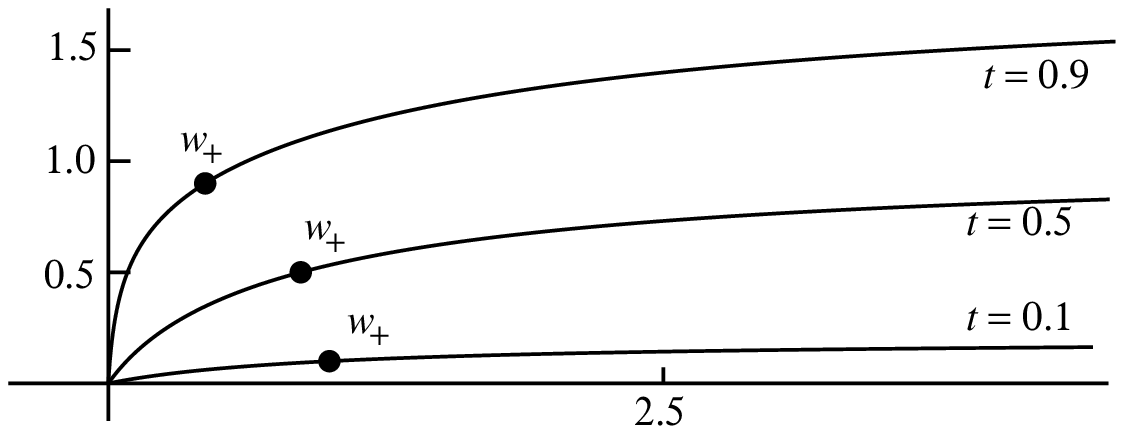}
\end{center}

\begin{quote}
{\bf Figure 2.}\quad
Steepest descent contours for the integral in \pref{nm:intw} for $t=0.1, 
0.5, 0.9$. 
\end{quote}

A simple approximation of the path is given by (we write $w=u+iv$)
\begin{equation}\label{nm:vofu}
v=\frac{ut(1+u_+)}{u+u_+^2},\quad u_+=\sqrt{1-t^2}.
\end{equation}
This path runs through the point $w_+=u_++it$, and has the same 
slope at this point as the exact steepest descent path, that is,  
$dv/du=t/(1+u_+)$ for $u=u_+$.
\vspace*{0.3cm}

For $t=1$ the steepest descent path runs from the origin to $w_+=i$ 
along the imaginary axis, 
and from $i$ to $2i+\infty$. For more details on the case $t\ge1$ we refer to 
\S\,\ref{subsec:uvnp}.

Integrating \pref{nm:intw} with respect to $\theta$
we obtain 
\begin{equation} \label{nm:intthuv}
\begin{array}{ll}
U(-a,x)+ & \dsp{i\Gamma (a+\tfrac12 )V(-a,x)=}\\
 &  \dsp{=\sqrt{\tfrac2\pi}\aof\gamma(a)
e^{-i(2a\eta-\frac14\pi)}} \,
\dsp{\int_0^{\theta_0} e^{-a\psi(\theta)} g(\theta)\,d\theta,}
\end{array}
\end{equation}
where $\gamma(a)$ is defined in \pref{pp:gammaa},
\begin{equation} \label{nm:psi}
\psi(\theta)=\tfrac12r^2\cos2\theta+2tr\sin\theta-\ln r-t^2-\tfrac12,
 \end{equation}
and
\begin{equation}\label{nm:gt}
g(\theta)=-\frac{e^{i\theta/2}}{\sqrt{r}}\left(\frac{dr}
{d\theta}+ir\right).  
 \end{equation}

We write the representations for $U(-a,x)$ and $V(-a,x)$ in real form, with
trigonometric functions that correspond with those in \cite{olpar} and
\cite{tempar}. We first write
\begin{equation}\label{nm:gt12}
g(\theta)=g_1(\theta)-ig_2(\theta),
 \end{equation}
where $g_j(\theta), j=1,2$, are real. That is, by \pref{nm:gt},
\begin{equation}\label{nm:gt1}
g_1(\theta)=-\frac{dr}{d\theta}\frac{\cos(\frac12\theta)}{\sqrt{r}}
+r\frac{\sin(\frac12\theta)}{\sqrt{r}},
 \end{equation}
\begin{equation}\label{nm:gt2}
g_2(\theta)=r\frac{\cos(\frac12\theta)}{\sqrt{r}}
+\frac{dr}{d\theta}\frac{\sin(\frac12\theta)}{\sqrt{r}}.
 \end{equation}
Then we have
\begin{equation}\label{nm:u12} 
U(-a,x)=
\sqrt{\tfrac2\pi}\aof\gamma(a)
\left[(\sin \lambda)\,G_1+(\cos\lambda )\,
G_2\right],
 \end{equation}
\begin{equation}\label{nm:v12} 
V(-a,x)=\sqrt{\tfrac2\pi}\aof
\frac{\gamma(a)}{\Gamma(a+\tfrac12)}
\left[(\cos\lambda )\,G_1-(\sin\lambda )
G_2\right],
 \end{equation}
where
\begin{equation}
\label{deflam}
\lambda=2a\eta+\tfrac14\pi
\end{equation}
\noindent and
\begin{equation} \label{nm:i12}
G_j=\int_{0}^{\theta_0} e^{-a\psi(\theta)}g_j(\theta)\,d\theta, \quad j=1,2.
 \end{equation}

For the derivatives we find, using \pref{nm:ints} and \pref{nm:yuv},
\begin{equation} \label{nm:intthuvd}
\begin{array}{ll}
&U'(-a,x)+i\Gamma\left(a+\tfrac12\right)V'(-a,x)=\\
{\ }\\
&\quad\quad
{\dsp
\sqrt{\tfrac2\pi}\atf\gamma(a)
e^{-i(2a\eta-\frac14\pi)} \,
\int_0^{\theta_0} e^{-a\psi(\theta)} h(\theta)\,d\theta.
}
\end{array}
\end{equation}
where
\begin{equation}\label{nm:ht}
h(\theta)=(t+iw)g(\theta)=h_1(\theta)-ih_2(\theta).  
 \end{equation}
That is, by \pref{nm:ht}, \pref{nm:gt1}
and \pref{nm:gt2},
\begin{equation}\label{nm:ht1}
h_1(\theta)=(t-r\sin\theta)g_1(\theta)+r\cos\theta g_2(\theta),
 \end{equation}
\begin{equation}\label{nm:ht2}
h_2(\theta)=r\cos\theta g_1(\theta)-(t-r\sin\theta) g_2(\theta).
 \end{equation}
Then we have
\begin{equation} \label{nm:u12j}
U'(-a,x)=
\sqrt{\tfrac2\pi}\atf\gamma(a)
\left[(\sin\lambda)\,H_1+(\cos\lambda)\,
H_2\right],
 \end{equation}
\begin{equation} \label{nm:v12j}
V'(-a,x)=
\sqrt{\tfrac2\pi}\frac{\atf \gamma(a)}{\Gamma(a+\frac12) }
\left[(\cos\lambda)\,H_1-(\sin\lambda)\,
H_2\right].
 \end{equation}
where $\lambda$ is given in Eq. (\ref{deflam}) and
\begin{equation} \label{nm:j12}
H_j=\int_{0}^{\pi} e^{-a\psi(\theta)}h_j(\theta)\,d\theta, \quad j=1,2.
 \end{equation}

\subsubsection{A Wronskian for the integrals}
\label{subsubsec:uvnm:wr}%
By using \pref{int:i10} and the integrals in \pref{nm:i12} and \pref{nm:j12}, 
we find the relation
\begin{equation} \label{nm:wij}  
H_1G_2-G_1H_2=\sqrt{\tfrac\pi2}\frac{\Gamma(a+\frac12)}
{a\gamma^2(a)}=\frac{\pi}{a}
 \Gamma^*\left(a+\tfrac12\right),
\end{equation}
where
\begin{equation}\label{nm:Gg} 
\Gamma(a+\tfrac12)=\sqrt{2\pi}\gamma^2(a)\Gamma^*(a+\tfrac12),\quad
\Gamma^*(a+\tfrac12)=1+\bo(1/a),
\end{equation}
as $a\to\infty$.
Hence, for large $a$, the right-hand side in \pref{nm:wij} is of order
$\frac{\pi}{a}[1+\bo(1/a)]$
(see also formula (3.28) in \cite{tempar}).
The relation in \pref{nm:wij} can be used for testing the numerical algorithms.

\subsubsection{Uniform asymptotic expansions {\protect\boldmath $-1<t<1$}}
\label{subsubsec:uvnm:as}%
The relationship of the integrals $G_j, H_j$ with uniform expansions 
follows from (2.23), (2.24) and (2.27) of \cite{tempar}. These expansions 
are the same as in \cite{olpar}. On the other hand, we can derive modified 
expansions (a main topic in \cite{tempar}), by using (2.29) and (2.33) of that 
reference. From \pref{nm:uuv} and by
changing $t$ to $-it$ in(2.29) of \cite{tempar},
we obtain
\begin{equation} \label{nm:yas}  
\begin{array}{ll}
&U(-a,x)+i\Gamma\left(a+\tfrac12\right)V(-a,x)=\\
{\ }\\
&\quad\quad
{\dsp
\Gamma^*\left(a+\tfrac12\right)\,
\frac{\sqrt{2}\gamma(a)e^{\frac14\pi 
i-2ia\eta}}{\aof(1-t^2)^{\frac14}}\,\wt{F}_{\mu}(-it),
}
\end{array}
\end{equation}
\begin{equation} \label{nm:ydas}  
\begin{array}{ll}
&U'(-a,x)+i\Gamma\left(a+\tfrac12\right)V'(-a,x)=\\
{\ }\\
&\quad\quad
{\dsp
-\Gamma^*\left(a+\tfrac12\right)\,
\sqrt{2}\aof\gamma(a)e^{-\frac14\pi 
i-2ia\eta}(1-t^2)^{\frac14}\,\wt{G}_{\mu}(-it),
}
\end{array}
\end{equation}
where $\wt{F}_{\mu}(-it)$ and $\wt{G}_{\mu}(-it)$ have the asymptotic 
expansions
\begin{equation} \label{nm:fas}  
\wt{F}_{\mu}(-it)\sim\sum_{s=0}^\infty (-1)^s\frac{\phi_s(\tau^*)}{(2a)^s},
\quad
\wt{G}_{\mu}(-it)\sim\sum_{s=0}^\infty (-1)^s\frac{\psi_s(\tau^*)}{(2a)^s},
\end{equation}
as  $a\to\infty$, uniformly for $t\in[-1+\delta, 1-\delta]$. The quantity 
$\tau^*$ is defined by 
\begin{equation} \label{nm:taustar}  
\tau^*=-\tfrac12\left(\frac{ it}{\sqrt{1-t^2}}+1\right).
\end{equation}
The polynomials $\phi_s$ and $\psi_s$ are 
given in (2.11) and (2.16) of \cite{tempar}, with recursion relations.
The first fraction at the right-hand sides of \pref{nm:yas} and
\pref{nm:ydas} has the 
asymptotic estimate $1+\bo(1/a)$ (see also formula (3.28) in \cite{tempar}).

\subsection{The case {\protect\boldmath $t\ge1$}}
\label{subsec:uvnp}%

We use the the integral for $Y(a,x)$ given in \pref{nm:intw}
with $\phi(w)$ given in \pref{nm:phiw}.
The saddle points are now purely imaginary:
\begin{equation}\label{np:wpm}
w_-= it-i\sqrt{t^2-1},\quad w_+=it+i\sqrt{t^2-1}.
\end{equation}
We have
\begin{equation}\label{np:phiwp}
\begin{array}{l}
\phi(w_\pm)=\tfrac12+t^2-\tfrac12\pi i \pm 2\xi, \\
\\
\xi = \tfrac12\left[t\sqrt{t^2-1}-\ln\left(t+ \sqrt{t^2-1}\right)\right].
\end{array}
 \end{equation}
The quantity $\xi$ is also is used in the asymptotic representation of $U(-a,x)$ 
for this case; 
see \cite{olpar} and \cite{tempar}.

The path of steepest descent starts at  $w=0$, runs through $w_-$ and $w_+$
on the positive imaginary axis,
and from $w_+$ to $+\infty$. 
The path from $w_+$ to $+\infty$ follows from solving the equation
\begin{equation}\label{np:imphi}
\Im\phi(w)=-\tfrac12\pi,
\end{equation}
that is, from solving
\begin{equation}\label{np:rt}
\tfrac12r^2\sin2\theta-2tr\cos\theta-\theta+\tfrac12\pi=0,
 \end{equation}
where $w=re^{i\theta}$.
The solution of \pref{np:rt} reads
\begin{equation}\label{np:r}
r=\frac{t\cos\theta+
\sqrt{t^2\cos^2\theta+\sin\theta\cos\theta(\theta-\frac12\pi)}}
{\sin\theta\cos\theta},\quad 0\le\theta\le\tfrac12\pi.
\end{equation}
The square root is positive, unless when $t=1$ and $\theta=\frac12\pi$. 

We obtain 
\begin{equation} \label{np:intthuv}
\begin{array}{ll}
U(-a,& x)+i\Gamma (a+1/2)V(-a,x)=\dsp{
\sqrt{\tfrac2\pi}\aof\gamma(a)\left[ e^{-2a\xi+\frac14\pi i} \right.}\\ 
\\
&\dsp{\times\left.
\int_0^{\frac12\pi} e^{-a\psi(\theta)} g(\theta)\,d\theta
+
ie^{2a\xi}\int_0^{r_+}e^{a\wt\phi(v)}\frac{dv}{\sqrt{v}}\right]},
\end{array}
\end{equation}
\noindent
where $\wt\phi(v)=\phi(w_-)-\phi(iv)$, 
$g(\theta)$ as in \pref{nm:gt},
and $\psi(\theta)=\phi(w)-\phi(w_+)=\Re[\phi(w)-\phi(w_+)]$, $w=r e^{i\theta}$,
now with $r$ defined in \pref{np:r}. Explicitly,
\begin{equation} \label{np:phivpsitheta}
\begin{array}{l}
\wt\phi(v)=\tfrac12v^2-2tv+\ln v- \tfrac12r_-^2+2tr_--\ln r_-,\\
\psi(\theta)=\tfrac12r^2\cos2\theta+2tr\sin\theta-\ln r-\tfrac12-t^2-2\xi.
\end{array}
 \end{equation}
\noindent
where
\begin{equation}
r_- =t-\sqrt{t^2-1}.
\end{equation}

Considering the real and imaginary parts on both sides of \pref{np:intthuv},
we see that 
for $V(-a,x)$ we need the $v-$integral with the dominant factor $e^{2a\xi}$
and part of the $\theta-$integral. When $t \sim 1$ (in fact, when $a\xi$ is 
small) both integrals
are of the same asymptotic importance. The dominant saddle point in the 
$v-$integral is $r_-$; in the $\theta-$integral the dominant point is 
the upper limit.

When we write
\begin{equation} \label{np:greim}
e^{\frac14\pi i}g(\theta)=g_1(\theta)+ig_2(\theta),
 \end{equation}
where $g_1(\theta)$ and $g_2(\theta)$ are real, we have
\begin{equation} \label{np:intthu}
U(-a,x)=
\sqrt{\tfrac2\pi}\aof\gamma(a)
e^{-2a\xi} \,G_1
\end{equation}
and
\begin{equation} \label{np:intthv}
V(-a,x)=\sqrt{\tfrac2\pi}\aof\frac{\gamma(a)e^{2a\xi}}{\Gamma\left(a+\frac12\right)}
\left(e^{-4a\xi} \,G_2+G_3\right).
\end{equation}
where (for $j=1,2$)
\begin{equation} \label{np:g123}
G_j=\int_0^{\frac12\pi} e^{-a\psi(\theta)} g_j(\theta)\,d\theta,
\quad
G_3=\int_0^{r_+}e^{a\wt\phi(v)}\frac{dv}{\sqrt{v}}.
\end{equation}

For the derivatives we have
\begin{equation} \label{np:intthuvd}
\begin{array}{ll}
U'(-a,& x) +  i\Gamma (a+\tfrac12 )  V'(-a,x)=\sqrt{\tfrac2\pi}\atf\gamma(a)\left[
e^{-2a\xi+\frac14\pi i} \right.\\
\\
& \dsp{\times\left.\int_0^{\frac12\pi} e^{-a\psi(\theta)} h(\theta)\,d\theta
+
ie^{2a\xi}\int_0^{r_+}e^{a\wt\phi(v)}(t-v)\frac{dv}{\sqrt{v}}\right]},
\end{array}
\end{equation}
where $h(\theta)=(t+iw)g(\theta)$. When we write
\begin{equation} \label{np:hreim}
e^{\frac14\pi i}h(\theta)=h_1(\theta)+ih_2(\theta),
 \end{equation}
where $h_1(\theta)$ and $h_2(\theta)$ are real, we have
\begin{equation} \label{np:intthud}
U'(-a,x)=
\sqrt{\tfrac2\pi}\atf\gamma(a)
e^{-2a\xi} \,H_1,
\end{equation}
and
\begin{equation} \label{np:intthvd} 
V'(-a,x)=
\sqrt{\tfrac2\pi}\atf\frac{\gamma(a)e^{2a\xi}}{\Gamma\left(a+\frac12\right)}
\left[e^{-4a\xi} \,H_2+H_3\right].
\end{equation}
where (for $j=1,2$)
\begin{equation} \label{np:h123}
H_j=\int_0^{\frac12\pi} e^{-a\psi(\theta)} h_j(\theta)\,d\theta,
\quad
H_3=\int_0^{r_+}e^{a\wt\phi(v)}(t-v)\frac{dv}{\sqrt{v}}.
\end{equation}

\subsubsection{A Wronskian for the integrals}
\label{subsubsec:uvnp:wr}%
By using the Wronskian relation in \pref{int:i10} and the integrals 
in \pref{np:g123} and  \pref{np:h123}, we obtain (cf. \pref{nm:wij})
\begin{equation}\label{np:wr1} 
\begin{array}{c}
e^{-4a\xi}(G_1H_2-H_1G_2)+(G_1H_3-H_1G_3)\\
\\
\dsp{=\sqrt{\tfrac\pi2}\frac{\Gamma(a+\frac12)}
{a\gamma^2(a)}=\frac{\pi}{a}\Gamma^* \left(a+\frac12\right)}.
\end{array}
\end{equation}
The relation in \pref{np:wr1} can be used for testing the numerical algorithms.

\subsubsection{Uniform asymptotic expansions  for {\protect\boldmath $t>1$}}
\label{subsubsec:uvnp:as}%
We give the relationship of the integrals with the uniform expansions 
given in (2.9), (2.14), (2.18) and (2.29)  of \cite{tempar}. We have
\begin{equation} \label{npn:asy} 
\begin{array}{l}
 G_1={\dsp\frac{\sqrt{\pi}}{2\sqrt{a} (t^2-1)^{\frac14}}\,{F}_\mu(t),}\\
\ \\
 H_1={\dsp-\frac{\sqrt{\pi}(t^2-1)^{\frac14}}{2\sqrt{a} }\,{G}_\mu(t),}\\
\ \\
 e^{-4a\xi} \,G_2+G_3={\dsp\frac{\Gamma(a+\tfrac12)e^aa^{-a-\frac12}}{\sqrt{2} 
(t^2-1)^{\frac14}}\,{P}_\mu(t),}\\
\ \\ 
e^{-4a\xi} \,H_2+H_3={\dsp\frac{\Gamma(a+\tfrac12)e^aa^{-a-\frac12}(t^2-1)^{\frac14}}{\sqrt{2} }
\,{Q}_\mu(t),}
\end{array}
\end{equation}
where ${F}_\mu(t)$, ${G}_\mu(t)$, ${P}_\mu(t)$, and ${Q}_\mu(t)$ are 
supplied with asymptotic expansions that have a double asymptotic property: 
one of the parameters $a$ or $t$ (or both) should be large; $t\ge1+\delta$. 
Recurrence 
relations for the coefficients of the expansions are given in \cite{tempar}.

\subsection{The case {\protect\boldmath $t\sim1$}}
\label{subsec:uvn1}%

For $t\sim 1$ the contours used in \S\,\ref{subsec:uvnm} becomes less
suitable for numerical quadrature.  For example, we see in Figure 2 that
the saddle point $w_+$ approaches the imaginary unit when $t\uparrow 1$,
and that the path becomes non-smooth when $t=1$.  For numerical
calculations we may consider uniform Airy-type asymptotic expansions if
$t\sim 1$, and we will investigate later if this is indeed the best
approach.  But we also investigate if a modified contour can be used for
numerical quadrature.

We use for $t\sim1$ the representation of $Y(a,x)$ in \pref{nm:intw}.
We write $w=u+iv$, and integrate with respect to $v$ along the line
segment from the origin to $w_+$, and then along the horizontal path
from $w_+$ to $w_++\infty$ with respect to $u$.  In the first
integral we substitute $v=t(1-p)$, and integrate with respect to $p$.
Observe that for $t\ge 1$ the point $w_+$ is on the imaginary axis,
and for this case no difficulties arise when $t\sim 1$, because the
path is already split up into two parts; see \S\,\ref{subsec:uvnp}.

It is not difficult to verify that the representations in \pref{nm:u12},
\pref{nm:v12}, \pref{nm:u12j} and \pref{nm:v12j} can be obtained, with $G_j, H_j$ 
replaced with $\wt G_j, \wt H_j$ ($j=1,2)$, where

\begin{equation}\label{n1:wtij}
\begin{array}{ll}
\wt G_j=\dsp{\int_0^1 e^{-a\psi_r^{(1)}(p)}} & \dsp{g_j^{(1)}(p)\,dp} \\&
\!\!\!\!\!\!\!\!\dsp{+
\int_0^\infty \frac{e^{-a\psi_r^{(2)}(u)}}
{(1+2u\sqrt{(1-t^2} +u^2)^\frac14}g_j^{(2)}(u)\,du,} \\ 
&
\\
\wt H_j=\dsp{\int_0^1 e^{-a\psi_r^{(1)}(p)}} & \dsp{h_j^{(1)}(p)\,dp}\\
&
\!\!\!\!\!\!\!\!\dsp{+
\int_0^\infty \frac{e^{-a\psi_r^{(2)}(u)}}
{(1+2u\sqrt{(1-t^2} +u^2)^\frac14}h_j^{(2)}(u)\,du,}
\end{array}
\end{equation}
where
\begin{equation}\label{n1:psiri}
\begin{array}{ll}
&\dsp{\psi_r^{(1)}(p)=\tfrac12p^2(1-2t^2)-p-\ln(1-p),}\\ \ \\
&\dsp{\psi_r^{(2)}(u)=\tfrac12u+u\sqrt{1-t^2}-
\tfrac12\ln\left(1+2u\sqrt{(1-t^2} +u^2\right),}\\ \ \\
&\dsp{\psi_i^{(1)}(p)=p^2 t\sqrt{1-t^2},}\\ \ \\
&\dsp{\psi_i^{(2)}(u)=\arctan\frac{ut}{1+u\sqrt{1-t^2}}-tu,}\\ \ \\
\end{array}
\end{equation}
\begin{equation}\label{n1:g}
\begin{array}{ll}
&\dsp{g^{(1)}_1(p)-ig^{(1)}_2(p)=e^{\frac12i\tau-ai\psi_i^{(1)}(p)},
\quad \tau=\arcsin t,}\\ \ \\
&\dsp{g^{(2)}_1(u)-ig^{(2)}_2(u)=e^{-\frac12i\tau-ai\psi_i^{(2)}(u)
+i\arctan\frac{ut}{1+u\sqrt{1-t^2}}},}\\ \ \\
\end{array}
\end{equation}
\begin{equation}\label{n1:h}
\begin{array}{ll}
&\dsp{h^{(1)}_1(p)=tpg_1^{(1)}(p)+\sqrt{1-t^2}(1-p)g_2^{(1)}(p),}\\ \ \\
&\dsp{
h^{(1)}_2(p)=tpg_2^{(1)}(p)-\sqrt{1-t^2}(1-p)g_1^{(1)}(p),}\\ \ \\
&\dsp{h^{(2)}_1(u)=(u+\sqrt{1-t^2})g_2^{(2)}(u),}\\ \ \\
&\dsp{
h^{(2)}_2(u)=-(u+\sqrt{1-t^2})g_1^{(2)}(u).}
\end{array}
\end{equation}

\subsection{The case {\protect\boldmath $t\le -1$}}
\label{subsec:uvnn}%

We can repeat the analysis, starting with \pref{nm:ints} with $x<0$,
but do not need new integral representations, algorithms or uniform 
asymptotic expansions 
for this case. For $U(a,x)$ we can use the 
second relation in \pref{int:i12}, and for $V(a,x)$ the first relation.

When the parameter $a$ is large these relations have to be used with care, 
because gamma functions with large negative arguments occur. It is better to use 
the quantities $G_j, H_j$ introduced in \S\,\ref{subsubsec:uvnp:wr}. 
In the computer code 
these quantities will be given as output from the case $t>1$.

We have
\begin{equation} \label{nn:mod} 
\begin{array}{l}
{\dsp U(-a,-x)=\sqrt{\tfrac2\pi}\aof\gamma(a) 
\left\{e^{-2a\xi}\left[\cos\pi a\, G_2 +\sin\pi a\, G_1\right]\right.} \\ 
\hspace*{6.3cm}{\dsp \left.+e^{2a\xi}\cos\pi a\, G_3\right\},}\\
\\
{\dsp U'(-a,-x)=-\sqrt{\tfrac2\pi}\atf\gamma(a)
\left\{e^{-2a\xi}\left[\cos\pi a\, H_2 +\sin\pi a\, H_1\right]\right.} \\
\hspace*{6.3cm}{\dsp \left.+e^{2a\xi}\cos\pi a\, H_3\right\},}\\
\\
{\dsp V(-a,-x)=\sqrt{\tfrac2\pi}\frac{\aof\gamma(a)}{\Gamma(a+\frac12)}
\left\{e^{-2a\xi}\left[\cos\pi a\, G_1 -\sin\pi a\, G_2\right]\right.}  \\
\hspace*{6.3cm}{\dsp \left. -e^{2a\xi}\sin\pi a\, G_3\right\},}\\
\\
{\dsp V'(-a,-x)=-\sqrt{\tfrac2\pi}\frac{\atf\,\gamma(a)}{\Gamma(a+\frac12)}
\left\{e^{-2a\xi}\left[\cos\pi a\, H_1 -\sin\pi a\, H_2\right]\right.}\\
\hspace*{6.3cm}{\dsp \left.-e^{2a\xi}\sin\pi a\, H_3\right\}}.
\end{array}
\end{equation}

\section{The {\protect\boldmath $W-$}function}
\label{sec:ww}
In this section solutions of equation
\begin{equation}\label{w:de}
\frac{d^2y}{dx^2}+\left(\tfrac14x^2-a\right)y=0.
\end{equation}
are considered, again for real $a$ and $x$. For $a<0$ the solutions 
oscillate on the real $x-$axis; for $a>0$ there are turning points at 
$\pm2\sqrt{a}$, and the oscillations occur outside the interval 
$[-2\sqrt{a},2\sqrt{a}]$. From quantum mechanics we know that 
\pref{w:de} is the equation for propagation through a potential barrier.

\subsection{The standard solutions}\label{sec:wpmx}
We consider solutions
$W(a,x)$ and $W(a,-x)$; these form a numerically
satisfactory pair for $-\infty<x<\infty$; see
\cite{mil55}.
The function $W(a,x)$ has the initial values (see \cite{abst}, p. 692)
\begin{equation}\label{w:init}
\begin{array}{l}
W(a,0)=2^{-\frac34}\left|\frac
{\Gamma(\tfrac14+\tfrac12ia)}
{\Gamma(\tfrac34+\tfrac12ia)}
\right|^{\frac12},
\\
\\
W'(a,0)=-2^{-\frac14}\left|\frac
{\Gamma(\tfrac34+\tfrac12ia)}
{\Gamma(\tfrac14+\tfrac12ia)}
\right|^{\frac12}.
\end{array}
\end{equation}
The Wronskian of $W(a,x)$ and $W(a,-x)$ is
 \begin{equation}\label{w:wr}
  {{\cal W}}[W(a,x),W(a,-x)] = 1.   
\end{equation}

Power series expansions are
\begin{equation}\label{w:ps}
W(a,x)=W(a,0)\,w_1(a,x)+W'(a,0)\,w_2(a,x), 
\end{equation}
where $w_1(a,x)$ and $w_2(a,x)$ are the
even and odd solutions of \pref{w:de}. We have
\begin{equation}\label{w:w12}
w_1(a,x)=\sum_{n=0}^\infty \alpha_n(a)\frac{x^{2n}}{(2n)!},
\quad
w_2(a,x)=\sum_{n=0}^\infty \beta_n(a)\frac{x^{2n+1}}{(2n+1)!},
\end{equation}
where $\alpha_n(a)$, $\beta_n(a)$ satisfy the recursion
\begin{equation}\label{w:al}
\begin{array}{ll}
&\alpha_{n+2} = a \,\alpha_{n+1} -\tfrac12 (n+1)(2n+1)\,\alpha_{n},\\
&\beta_{n+2} = a \,\beta_{n+1} -\tfrac12 (n+1)(2n+3)\,\beta_{n}, \\
&\alpha_0(a) = 1,\quad \alpha_1(a)=a,\quad \beta_0(a)=1, \quad
\beta_1(a)=a.
\end{array}
\end{equation} 

The relation with the function $U(a,x)$ reads
\begin{equation}\label{w:uu}
\begin{array}{lr}
\frac{1}{\sqrt{k(a)}}W(a,x)+i\sqrt{k(a)} W(a,-x)&\\
& \hspace*{-2cm}=\sqrt{2}
e^{\frac14\pi a+i\rho(a)}
U\left(ia,xe^{-\pi i/4}\right),
\end{array}
\end{equation}
which follows from using the initial values of the functions, but also from
\cite{abst} [19.17.6 and 19.17.9]. The quantities $k(a)$ and $\rho(a)$ are 
given by
\begin{equation}\label{w:k}
k(a)=\sqrt{1+e^{2\pi a}}-e^{\pi a} 
= \frac{1}{\sqrt{1+e^{2\pi a}}+e^{\pi a}},
\end{equation}
and
\begin{equation}\label{w:rho}
\rho(a)=\tfrac18\pi+\tfrac12\phi_2(a),
\quad
\phi_2(a)=\ph\,{\Gamma\left(\tfrac12+ia\right)};
\end{equation}
the branch is defined  by $\phi_2(0)=0$ and by
continuity elsewhere. 

Because we assume that $a$ and $x$, and hence $W(a,\pm x)$, are real,
we have, using \pref{w:uu}, that
\begin{equation}\label{w:pmx}
\begin{array}{ll}
&\dsp{W(a,x)=
\sqrt{2k(a)}
e^{\frac14\pi a} \Re\left[ e^{i\rho(a)}
U\left(ia,xe^{-\pi i/4}\right)\right],}\\
&\  \\
&\dsp{W(a,-x)=
\sqrt{\frac{2}{k(a)}}
e^{\frac14\pi a} \Im\left[ e^{i\rho(a)}
U\left(ia,xe^{-\pi i/4}\right)\right].}
\end{array}
\end{equation}
These relations are convenient for numerical computations because for $x\ge0$ 
and $x\le0$ we can use the same $U-$function.

\subsubsection{The function  {\protect\boldmath $\rho(a)$} }
We give more details on the function $\rho$ defined in \pref{w:rho}.
For large values of $a$ it is convenient to use the 
representation
\begin{equation}\label{w:rhostar}
\rho(a)=\tfrac18\pi-\tfrac12a+\tfrac14a\ln a^2+\rho^*(a),
\end{equation}
where $\rho^*(a)=\bo(1/a)$ as $a\to\infty$. 
To give more details we give an 
asymptotic expansion. We have Binet's formula (see \cite{temsf}, p. 55, 
for an integrated version)
\begin{equation}\label{w:binet}
\begin{array}{ll}
\ln\Gamma(z+\tfrac12)= & z\ln(\tfrac12+z)-\tfrac12-z+\tfrac12\ln(2\pi)\\
&\\
 & {\dsp +\int_0^\infty \beta(t) e^{-\frac12t} \, e^{-zt}\,dt},
\end{array}
\end{equation}
where
\begin{equation}\label{w:beta}
\begin{array}{ll}
\beta(t)\,e^{-\frac12t}&={\dsp t^{-1}\left(\frac1{e^t-1}-\frac1t+
\frac12\right) 
\,e^{-\frac12t}} \\
&\\
& {\dsp =\sum_{k=0}^\infty 
\frac{c_k\,t^k}{(k+2)!},\quad |t|<2\pi},
\end{array}
\end{equation}
with $c_k$ in terms of Bernoulli polynomials:
\begin{equation}\label{w:ck}
c_k=B_{k+2}^{(1)}(-\tfrac12)-(-1)^k\frac{k+3}{2^{k+2}}, \quad k=0,1,2,\ldots\ .
\end{equation}
This gives the asymptotic expansion
\begin{equation}\label{w:asrho}
\rho^*(a)\sim\tfrac14a\ln\left(1+\frac1{4a^2}\right)-\frac1{2a}\sum_{k=0}^\infty
\frac{d_k}{a^{2k}},
\end{equation}
as $\pm a\to\infty$, where
\begin{equation}\label{w:dk}
d_k=(-1)^k\frac{c_{2k}}
{(2k+1)(2k+2)}, \quad k=0,1,2,\ldots\ . 
\end{equation}\
The first few coefficients are
\begin{equation}\label{w:dk04}
\hspace*{-0.1cm}d_0= \tfrac{1}{12},\ 
d_1= -\tfrac{13}{720},\ 
d_2= \tfrac{37}{20160},\ 
d_3= -\tfrac{29}{26880},\ 
d_4= -\tfrac{1129}{1520640}.
 \end{equation}\

\subsection{Integral representations for {\protect\boldmath $a<0$} }
\label{subsec:wwn}%
For $W(-a,\pm x)$ we consider 
\pref{pp:int} for 
$U\left(-ia,xe^{-\pi  i/4}\right)$ (see \pref{w:uu} and \pref{w:pmx}), 
that is,
\begin{equation}\label{wnp:ints}
U\left(-ia,xe^{-\pi  i/4}\right)=
\frac{e^{-\frac14ix^2}}{i\sqrt{2\pi}}\int_{\C} e^{-xe^{-\pi  i/4}s+
\frac12s^2}s^{ia}
\frac{ds}{\sqrt{s}},
\end{equation}
where $\C$ is a vertical line on which $\Re{s}>0$. 
On $\C$ we have $-\frac12\pi<\ph{\,s}<\frac12\pi$,
and the many-valued function $s^{ia-\frac12}$ assumes its principal value.
The transformations
\begin{equation}\label{wnp:tr}
x=2t\sqrt{a},\quad s=\sqrt{a}\,w
\end{equation}
give
\begin{equation}\label{wnp:intw}
U\left(-ia,xe^{-\pi  i/4}\right)=
\frac{e^{-\frac14ix^2}a^{\frac14+\frac12 ai}}{i\sqrt{2\pi}}\int_{\C} 
e^{a\phi(w)}
\frac{dw}{\sqrt{w}},
\end{equation}
where
\begin{equation}\label{wnp:phi}
\phi(w)=\tfrac12w^2-2te^{-\pi  i/4}w+i\ln w.
\end{equation}
The saddle points follow from solving
\begin{equation}\label{wnp:phip}
\phi'(w)=\frac{w^2-2te^{-\pi  i/4}w+i}{w}=0,
\end{equation}
giving saddle points
\begin{equation}\label{wnp:wpm}
w_{\pm}=u_{\pm}+iv_{\pm}=e^{-\pi  i/4}\left(t\pm\sqrt{t^2+1}\right), \quad 
v_{\pm}=-u_{\pm}. 
\end{equation}
The relevant saddle point is $w_+$. We have
\begin{equation}\label{wnp:phiwp}
\phi(w_+)=\tfrac14\pi+i(t^2 +2\wt\xi-\tfrac12),
\end{equation}
where $\wt\xi$ is given in \pref{pp:xit}. The path of steepest descent 
through $w_+$ is for $|\theta|<\tfrac12\pi$ defined by
\begin{equation}\label{wnp:rt} 
\tfrac12 r^2\sin2\theta -2tr\sin(\theta-\tfrac14\pi)+\ln r 
=\Im\phi(w_+)=t^2 +2\wt\xi-\tfrac12,
\end{equation} 
where $w=re^{i\theta}$. 
In rectangular coordinates $w=u+iv$ this equation  reads
\begin{equation}\label{wnp:uv} 
uv +\sqrt{2}t(u-v)+\tfrac12\ln(u^2+v^2) = t^2 +2\wt\xi-\tfrac12.
\end{equation} 
We can solve equation \pref{wnp:rt} for $\sin (\theta-\frac14\pi)$ 
(it is a quadratic equation for this quantity), giving $\theta$ as function of 
$r$.
This makes it possible to integrate \pref{wnp:intw} with respect to $r$,
but this introduces singularities in the integral where $r$ attains 
its minimal value, although the path itself is smooth.

Integrating with respect to $\theta$ or $v$ is a better option.
We can numerically determine the path in an algorithm, but this is not a very 
efficient method.  
Instead,  we replace the steepest 
descent path defined in \pref{wnp:uv} by a path $u(v)$ such that  
\begin{enumerate}
\item  $u(v)$ is smooth for all $v\in\RR$;
\item $u(v)$ passes through the saddle point: $u(v_+)=u_+$;
\item $du/dv$ at $v_+$ has the same value as $du/dv$ 
for the steepest descent contour at $w_+$;
\item the path runs into the valleys of $e^{\phi(w)}$ at 
                $\pm i\infty$.
\end{enumerate}
\noindent
>From \pref{wnp:uv} we can show that $du/dv=0$  at the saddle point $w_+$.
Hence, a simple path $\C$ that fulfills the four conditions is 
the vertical line $u=u_+$.
Introducing $q=v-v_+$, using $w=w_+ +iq$ and 
\pref{wnp:phiwp} we obtain for  \pref{wnp:intw}  the representation
\begin{equation}\label{wnp:intq}
\begin{array}{ll}
U\left(-ia,xe^{-\pi  i/4}\right)=&{\dsp
\frac{e^{\frac14\pi a-\frac12ia}a^{\frac14+\frac12 
ai}}{\sqrt{2\pi w_+}} \,e^{2ia\wt\xi}\, }\\
& \\
& \hspace*{0.5cm}{\dsp \times \int_{-\infty} ^\infty
e^{-a\psi(q)} g(q)\,dq.}
\end{array}
\end{equation}
The function $\psi(q)$ is given by
\begin{equation}\label{wnp:psiq}
\psi(q)=\phi(w_+)-\phi(w), \quad g(q)= 
\frac{1}{\sqrt{1+ iq/w_+}}.
\end{equation}
For small values of $q$ we have
\begin{equation}\label{wnp:psiexp}
\begin{array}{l}
\psi(q)=\frac{1+2u_+^2}{4u_+^2}q^2+
\frac{1-i}{12u_+^3}q^3+\bo(q^4),\\
\\
u_+=\frac1{\sqrt{2}}\left(t+\sqrt{t^2+1}\right).
\end{array}
\end{equation}

We conclude that  $W(-a,\pm x)$ (see \pref{w:uu} -- \pref{w:rhostar}) 
are given by
\begin{equation}\label{wnp:re}
W(-a,x)=\frac{a^{\frac14}\sqrt{k(-a)}}{\sqrt{\pi\,|w_+|}}\,
\Re\left[e^{i\chi}\,
\int_{-\infty} ^\infty e^{-a\psi(q)} g(q)\,dq\right],
\end{equation}
\begin{equation}\label{wnp:im}
W(-a,-x)=\frac{a^{\frac14}}{\sqrt{\pi\,|w_+|k(-a)}}\,
\Im\left[e^{i\chi}\,
\int_{-\infty} ^\infty e^{-a\psi(q)} g(q)\,dq\right],
\end{equation}
where
\begin{equation}\label{wnp:chi}
\chi=\rho^*(-a)+\frac14\pi+2a\wt\xi.
\end{equation}
For the derivatives we find, starting with \pref{wnp:ints},
\begin{equation}\label{wnp:intdq}
\begin{array}{ll}
e^{-\pi  i/4}U'\left(-ia,xe^{-\pi  i/4}\right)=&
{\dsp i\frac{e^{\frac14\pi a-\frac12ia}a^{\frac34+\frac12 
ai}}{\sqrt{2\pi w_+}} \,e^{2ia\wt\xi}\,} \\
\\
&\hspace*{1cm}{\dsp \times \int_{-\infty} ^\infty
e^{-a\psi(q)} h(q)\,dq,}
\end{array}
\end{equation}
where
\begin{equation}\label{wnp:gq}
h(q) =\left(\sqrt{t^2+1}-e^{-\pi i/4}q\right) g(q).
\end{equation}
It follows from \pref{w:pmx} that  $W'(-a,\pm x)$ are given by
\begin{equation}\label{wnp:red}
\hspace*{-1.2cm}W'(-a,x)=\frac{a^{\frac34}\sqrt{k(-a)}}{\sqrt{\pi\,|w_+|}}\,
\Re\left[i\,e^{i\chi}\,
\int_{-\infty} ^\infty e^{-a\psi(q)} h(q)\,dq\right],
\end{equation}
\begin{equation}\label{wnp:imd}
\hspace*{-0.1cm}W'(-a,-x)=-\frac{a^{\frac34}}{\sqrt{\pi\,|w_+|k(-a)}}\,
\Im\left[i\,e^{i\chi}\,
\int_{-\infty} ^\infty e^{-a\psi(q)} h(q)\,dq\right].
\end{equation}
 
For large values of $a$ and/or $t$ the oscillatory behaviour of $W(-a,\pm x)$ and 
$W'(-a, \pm x)$ is mainly described by the exponential factor $e^{2ia\wt\xi}$
contained in $e^{i\chi}$. The other elements of these formulas 
are slowly varying. 

Asymptotic expansions follow from \cite{tempar}. Formula (2.29) 
of that paper gives, with $\mu=\sqrt{2a}\,e^{-\pi i/4}$, 
\begin{equation}\label{wnp:was}
W(-a,x)\sim 
\frac{\sqrt{k(-a)}}{a^{\frac14}(t^2+1)^{\frac14}}\,
\Re\left[e^{i\chi}
\,
\sum_{s=0}^\infty\frac{(-i)^s\phi_s(\wt\tau)}{(2a)^s}\right],
\end{equation}
where $\phi_s$ are polynomials given in (2.11)  and 
$\wt\tau$ in (2.32) of \cite{tempar}.

Formula (2.33) of \cite{tempar} gives 
\begin{equation}\label{wnp:wasd}
W'(-a,x)\sim 
\sqrt{k(-a)}\,a^{\frac14}(t^2+1)^{\frac14}\,
\Re\left[i\,e^{i\chi}\,
\sum_{s=0}^\infty\frac{(-i)^s\psi_s(\wt\tau)}{(2a)^s}\right],
\end{equation}
where  $\psi_s$ are polynomials given in (2.16) of \cite{tempar}.

For $W(-a,-x)$ and its derivative we have
\begin{equation}\label{wnn:was}
W(-a,-x)\sim 
\frac{\sqrt{k(-a)}}{a^{\frac14}(t^2+1)^{\frac14}}\,
\Im\left[e^{i\chi}
\,
\sum_{s=0}^\infty\frac{(-i)^s\phi_s(\wt\tau)}{(2a)^s}\right],
\end{equation}
and
\begin{equation}\label{wnn:wasd}
\hspace*{-0.2cm} W'(-a,-x)\sim 
-\sqrt{k(-a)}\,a^{\frac14}(t^2+1)^{\frac14}\,
\Im\left[i\,e^{i\chi}\,
\sum_{s=0}^\infty\frac{(-i)^s\psi_s(\wt\tau)}{(2a)^s}\right],
\end{equation}

The asymptotic expansions in \pref{wnp:was} -- \pref{wnn:wasd} hold when 
$a\to\infty$, uniformly with respect to $t\ge -t_0$, but also for $t\to\infty$, 
uniformly 
with respect to $a\ge a_0$, where $a_0$ and $t_0$ are fixed positive 
numbers. 

\subsection{Integral representations for {\protect\boldmath $a>0$} }
\label{subsec:wwp}%

Because of the turning points we consider three cases. We write $x=2t\sqrt{a}$.
We use the  $U-$function in \pref{w:uu}, and write \pref{pp:int}
in the form
\begin{equation}\label{wp:ints}
U\left(ia,xe^{-\pi  i/4}\right)=
\frac{e^{-\frac14ix^2}}{i\sqrt{2\pi}}\int_{\C} e^{-xe^{-\pi  i/4}s+
\frac12s^2}s^{-ia}
\frac{ds}{\sqrt{s}},
\end{equation}
with conditions as in \pref{wnp:ints}. The transformation $s=\sqrt{a}w$ gives
\begin{equation}\label{wp:intw}
U\left(ia,xe^{-\pi  i/4}\right)=
\frac{e^{-\frac14ix^2}a^{\frac14-\frac12 ai}}{i\sqrt{2\pi}}\int_{\C} 
e^{a\phi(w)}
\frac{dw}{\sqrt{w}},
\end{equation}
where
\begin{equation}\label{wp:phi}
\phi(w)=\tfrac12w^2-2te^{-\pi  i/4}w-i\ln w.
\end{equation}

\subsubsection{The case {\protect\boldmath $t\ge1$} }
\label{subsubsec:wwpp}%
The saddle points are now
\begin{equation}\label{wpp:wpm}
w_{\pm}=u_{\pm}+iv_{\pm}=e^{-\pi  i/4}\left(t\pm\sqrt{t^2-1}\right), \quad 
v_{\pm}=-u_{\pm}. 
\end{equation}
The relevant saddle point is $w_+$, and for numerical integration a 
convenient choice of $\C$ is the vertical line through $w_+$.

Using
\begin{equation}\label{wpp:xi}
-\tfrac14ix^2+a\phi(w_+)=2ia\xi+\tfrac12ia-\tfrac14\pi a,
\end{equation}
where $\xi$ is given in \pref{np:phiwp}, and  writing in \pref{wp:intw} 
$w=w_++iq$, we obtain
the analogue of \pref{wnp:intq}
\begin{equation}\label{wp:intq}
U\left(ia,xe^{-\pi  i/4}\right)=
\frac{e^{-\frac14\pi a+\frac12ia}a^{\frac14-\frac12 
ai}}{\sqrt{2\pi w_+}} \,e^{2ia\xi}\, \int_{-\infty} ^\infty
e^{-a\psi(q)} g(q)\,dq,
\end{equation}
where
\begin{equation}\label{wpp:psiq}
\psi(q)=\phi(w_+)-\phi(w), \quad g(q)= 
\frac{1}{\sqrt{1+ iq/w_+}}.
\end{equation}
It follows that 
\begin{equation}\label{wpp:re}
W(a,x)=\frac{\sqrt{k(a)}\,a^{\frac14}}{\sqrt{\pi\,|w_+|}}\,
\Re\left[ F(a,x)\right] 
\end{equation}
and
\begin{equation}\label{wpp:im}
W(a,-x)=\frac{\sqrt{k(a)}\,a^{\frac14}}{\sqrt{\pi\,|w_+|}}\,
\Im\left[ F(a,x)\right].
\end{equation}
\noindent
where 
\begin{equation}
F(a,x)= e^{i[\rho^*(a)+\frac14\pi+2ia\xi]}
 \int_{-\infty} ^\infty e^{-a\psi(q)} g(q)\,dq
\end{equation}

For the derivative we find, as in \pref{wnp:intdq},
\begin{equation}\label{wpp:intdq}
\begin{array}{ll}
U'\left(ia,xe^{-\pi  i/4}\right)=&{\dsp
-\frac{e^{-\frac14\pi a+\frac12ia}a^{\frac34-\frac12 
ai}}{\sqrt{2\pi w_+}} \,e^{2ia\xi}}\,\\
\\
&{\dsp \hspace*{1cm}\times \int_{-\infty} ^\infty
e^{-a\psi(q)} h(q)\,dq,}
\end{array}
\end{equation}
where
\begin{equation}\label{wpp:gq}
h(q) =\left(e^{-\pi i/4}\sqrt{t^2-1}+iq\right) g(q).
\end{equation}
It follows that 
\begin{equation}\label{wpp:red}
\hspace*{-0.3cm}W'(a,x)=-\frac{\sqrt{k(a)}\,a^{\frac34}}{\sqrt{\pi}|w_+|}\,
\Re\left[ G(a,x)\right]
\end{equation}
and
\begin{equation}\label{wpp:imd}
W'(a,-x)=\frac{a^{\frac34}}{\sqrt{\pi k(a)}\,|w_+|}\,
\Im\left[ G(a,x)\right].
\end{equation}
\noindent
where
\begin{equation}
G(a,x)=e^{i[\rho^*(a)+\frac14\pi+2ia\xi]}\,
 \int_{-\infty} ^\infty e^{-a\psi(q)} h(q)\,dq
\end{equation}

Asymptotic expansions follow from (2.9) of \cite{tempar}. By changing
$\mu\to \mu e^{-\pi i/4}$ in that formula we obtain
\begin{equation}\label{wpp:was}
W(a,x)\sim 
\frac{\sqrt{k(a)}}
{a^{\frac14}(t^2-1)^{\frac14}}\,
\Re\left[e^{i[\rho^*(a)+\frac14\pi+2a\xi]}\,
\sum_{s=0}^\infty\frac{i^s\phi_s(\tau)}{(2a)^s}\right],
\end{equation}
where $\phi_s(\tau)$ are the same polynomials as in \pref{wnp:was}, and
\begin{equation}\label{wpp:tau}
\tau=\frac{1}{2}\left(\frac{t}{\sqrt{t^2-1}}-1\right).
\end{equation}

Formula (2.18) of  \cite{tempar} gives 
\begin{equation}\label{wpp:wasd}
\begin{array}{ll}
W'(a,x)\sim &
-\sqrt{k(a)}\,a^{\frac14}(t^2-1)^{\frac14}\\
\\
& {\dsp \vspace*{3cm}\times\Re\left[
e^{i[\rho^*(a)-\frac14\pi+2a\xi]}\,
\sum_{s=0}^\infty\frac{i^s\psi_s(\tau)}{(2a)^s}\right]},
\end{array}
\end{equation}
where $\psi_s(\tau)$ are the same as in \pref{wnp:wasd}.

The asymptotic expansions in \pref{wpp:was} and \pref{wpp:wasd} hold when 
$a\to\infty$, uniformly with respect to $t\ge 1+t_0$, but also for $t\to\infty$, 
uniformly 
with respect to $a\ge a_0$, where $a_0$ and $t_0$ are   positive 
numbers. 

\subsubsection{The case {\protect\boldmath $-1\le t\le1$} }
\label{subsubsec:wwpm}%
We use \pref{wp:ints}, \pref{wp:intw} and \pref{wp:phi} with saddle 
points
\begin{equation}\label{wpm:wpm}
w_{\pm}=e^{-\pi  i/4}\left(t\pm i\sqrt{1-t^2}\right) = 
e^{-\pi  i/4\pm i\theta}, \quad t=\cos\theta,
\end{equation}
which are located on the unit circle. We have
\begin{equation}\label{wpm:phis}
\phi(w_\pm)= \pm 2\eta - \tfrac14\pi+i(\tfrac12+t^2),
\end{equation}
where $\eta=\frac12(\theta-\sin\theta\cos\theta)$ is also used in 
\S \ref{subsec:uvnm}
and defined in \pref{nm:phiwp}.  
We see that the imaginary parts of 
$\phi(w_\pm)$ are equal. As a consequence, the steepest descent path may go
(and in fact in the present case does go)  
through both saddle points.

In Figure 3 we have shown the paths for three 
values of $t$.  The contours run from $-i\infty$ to $w_-$, 
then along the arc to $w_+$ (in the 
direction of the arrows), and from $w_+$ to $+i\infty$. 
Through each saddle point the local contours of 
steepest ascent and steepest ascent are shown. The complete contours include
steepest descent parts and steepest ascent parts.

>From \pref{wpm:phis} we see that $w_+$ is dominant for $0\le t<1$ ($\eta$ is 
positive for these values of $t$). Another point of interest is that the 
oscillatory factor $e^{-\frac14ix^2}=e^{-iat^2}$ in \pref{wp:intw} is nullified 
when we put $\phi(w_-)$ or $\phi(w_+)$ in front of the integral, 
because $\Im\phi(w_\pm)=\frac12+t^2$. This explains that the function 
$W(a,x)$ does not oscillate  if $t\in[-1,1]$.

When we write $w=u+iv$ and integrate in \pref{wp:intw} with respect to $v$
we obtain
\begin{equation}\label{wpm:intv}
\begin{array}{ll}
&\dsp{U\left(ia,xe^{-\pi  i/4}\right)=
\frac{a^{\frac14-\frac12 ai}e^{a[\frac12i-\frac14\pi+2\eta]}}
{\sqrt{2\pi}}\ \times\ }\\
&\  \\
&\quad
{\dsp{\left[\int_{v_-}^\infty  e^{-a\psi_1(v)} f_1(v)\,dv
+e^{-4a\eta}\int_{-\infty} ^{v_-} e^{-a\psi_2(v)} f_2(v)\,dv\right],}}
\end{array}
\end{equation}
where
\begin{equation}\label{wpm:psi}
\psi_1(v)= \phi(w_+)-\phi(w), \quad \psi_2(v)=\phi(w_-)-\phi(w),
\end{equation}
and
\begin{equation}\label{wpm:fv}
f_1(v)=\frac1{\sqrt{w}}\left(1-i\frac{du_1}{dv}\right), 
\quad f_2(v)=\frac1{\sqrt{w}}\left(1-i\frac{du_2}{dv}\right).
\end{equation}
We may assume different relations between $u$ and $v$ in both integrals; this
explains $u_1$ and $u_2$, which are functions of $v$.

\begin{center}
\epsfxsize=11cm \epsfbox{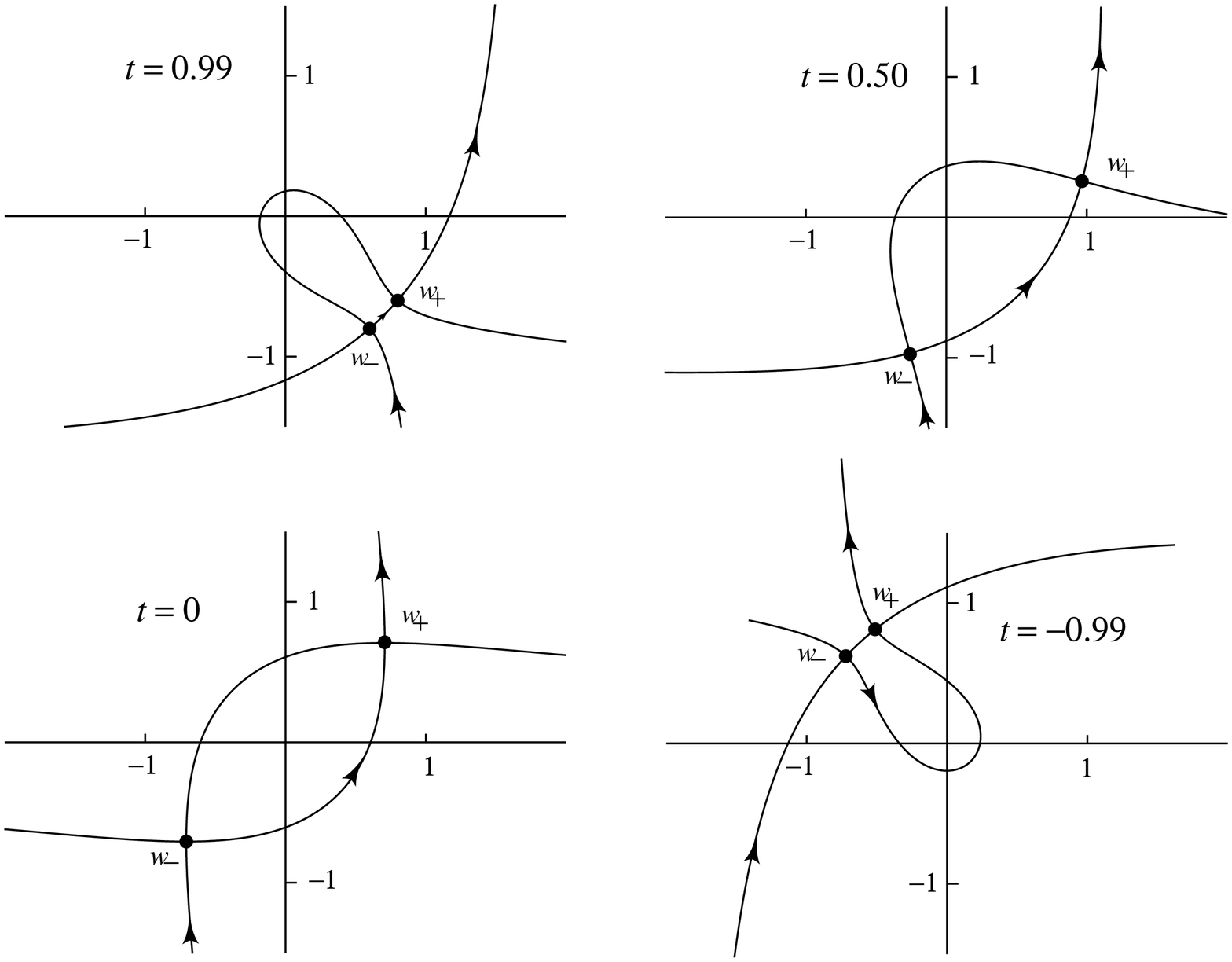}
\end{center}

\begin{quote}
{\bf Figure 3.}\quad
{
Steepest descent contour for the integral in \pref{wpm:intv} 
for several $t-$values. The contours run from $-i\infty$ to $w_-$, 
then along the arc to $w_+$ (in the direction of the arrow), 
and from $w_+$ to $+i\infty$. Through each saddle point the local contours of 
steepest ascent and steepest ascent are shown. The complete contours include
steepest descent parts and steepest ascent parts.
}
\end{quote}

It follows from \pref{w:pmx} that
\begin{equation}\label{wpm:reim}
\begin{array}{ll}
&\dsp{W(a,x)=\frac{\sqrt{k(a)}\,a^{\frac14}e^{2a\eta}}{\sqrt{\pi}}\,
\Re\left[ e^{i[\rho^*(a)+\frac18\pi]}\,K(a,x)\right],}\\
&\  \\
&\dsp{W(a,-x)=\frac{a^{\frac14}e^{2a\eta}}{\sqrt{\pi\,k(a)}}\,
\Im\left[ e^{i[\rho^*(a)+\frac18\pi]}\,K(a,x)\right],}
\end{array}
\end{equation}
where $K(a,x)$ denotes  the sum of the integrals between the square 
brackets in \pref{wpm:intv}.

For the derivatives we find, starting with \pref{wp:intw},
\begin{equation}\label{wpm:intdw}
\hspace*{-0.2cm}U'\left(ia,xe^{-\pi  i/4}\right)=
\frac{e^{-\frac14ix^2}a^{\frac34-\frac12 ai}}{\sqrt{2\pi}}\int_{\C} 
e^{a\phi(w)}\left(te^{-\pi  i/4}-w\right)
\frac{dw}{\sqrt{w}},
\end{equation}
and using \pref{w:pmx},
\begin{equation}\label{wpm:red}
\begin{array}{ll}
&\dsp{W'(a,x)=\frac{\sqrt{k(a)}\,a^{\frac34}e^{2a\eta}}{\sqrt{\pi}}\,
\Re\left[ e^{i[\rho^*(a)-\frac18\pi]}\,
 K_d(a,x)\right],}\\
&\ \\
&\dsp{W'(a,-x)=-\frac{a^{\frac34}e^{2a\eta}}{\sqrt{\pi k(a)}}\,
\Im\left[ e^{i[\rho^*(a)-\frac18\pi]}\,
 K_d(a,x)\right],}
\end{array}
\end{equation}
where $K_d(a,x)$ denotes  the sum of the integrals between the square 
brackets in \pref{wpm:intv} with $f_i(v)$ replaced with 
$g_i(v)=(te^{-\frac14\pi i}-w)f_i(v)$.

The asymptotic expansions follows from (2.29) and (2.33) of \cite{tempar}. 
We change 
$t\to-it$ and $\mu \to \mu e^{\frac14\pi i}$ in these formulae and obtain
\begin{equation}\label{wpm:uas}
\begin{array}{l}
\dsp{U\left(ia,xe^{-\pi  
i/4}\right)\sim
\frac{a^{-\frac12ia-\frac14}e^{\frac12ia-\frac18\pi i-\frac14\pi a+2a\eta}}
{\sqrt{2}(1-t^2)^{\frac14}} \sum_{s=0}^\infty 
i^s\frac{\phi_s(\tau^*)}{(2a)^s},}\\
 \\
\dsp{U'\left(ia,xe^{-\pi  
i/4}\right)\sim
-\frac{a^{-\frac12ia+\frac14}e^{\frac12ia+\frac18\pi i-\frac14\pi a+2a\eta}
(1-t^2)^{\frac14}}
{\sqrt{2}}}\\
\hspace*{5cm}{\dsp \times\sum_{s=0}^\infty 
i^s\frac{\psi_s(\tau^*)}{(2a)^s},}
\end{array}
\end{equation}
where $\eta$ is defined in \pref{nm:phiwp} and $\phi_s(\tau^*)$ are the 
same as in \pref{nm:fas}. 

It follows from \pref{w:pmx} that
\begin{equation}\label{wpm:was}
\begin{array}{l}
\dsp{W(a,x)\sim 
\frac{\sqrt{k(a)}}
{a^{\frac14}(1-t^2)^{\frac14}}e^{2a\eta}\,
\Re\left[e^{i\rho^*(a)}\,
\sum_{s=0}^\infty\frac{i^s\phi_s(\tau^*)}{(2a)^s}\right],}\\
 \\
\dsp{W'(a,x)\sim 
-\sqrt{k(a)}\,a^{\frac14}(1-t^2)^{\frac14}e^{2a\eta}}\\
\hspace*{3cm}\times{\dsp\Re\left[
e^{i[\rho^*(a)+\frac14\pi]}\,
\sum_{s=0}^\infty\frac{i^s\psi_s(\tau^*)}{(2a)^s}\right].}
\end{array}
\end{equation}
The asymptotic expansions in \pref{wpm:was}  hold when 
$a\to\infty$, uniformly with respect to $t\in[-1+\delta,1-\delta]$, where
$\delta$ is a fixed positive number. 

\subsubsection{Unstable representations}
\label{subsubsec:unst}%
For large values of $a\eta$ the representations for $W(a,-x)$ and $W'(a,-x)$ in 
\pref{wpm:reim} and \pref{wpm:red} are unstable. To see this, observe that 
\pref{wpm:was} can be used for $t\in[-1+\delta,1-\delta]$. 
The dominant behaviour comes from $\sqrt{k(a)} e^{2a\eta}$. Since (see \pref{w:k})
$k(a)\sim \frac12e^{-\pi a}$, the dominant behaviour comes from 
$e^{a\chi}$, where $\chi=\arcsin t-t\sqrt{1-t^2}$, an odd function that 
is positive on $(0,1]$. This dominant behaviour does not appear in
the representations for $W(a,-x)$ and $W'(a,-x)$ in 
\pref{wpm:reim} and \pref{wpm:red}. There we see the dominant parts
$ e^{2a\eta}/\sqrt{k(a)}$; in $\eta$ we use positive $t$ when $x$ is positive. 
It follows that the imaginary parts in the right-hand sides of 
\pref{wpm:reim} and \pref{wpm:red} have to be very small when $a\eta$ is large. 
In fact, the first integral in \pref{wpm:intv} should be of order $e^{-4a\eta}$
in that case, which is not apparent from this representation.
 
A possible solution to this problem is using
the representations for $W(a,x)$ and $W'(a,x)$ in 
\pref{wpm:reim} and \pref{wpm:red} for  $t\in[-1,0]$. However, when 
$t\downarrow -1$ the phase of $w_+$ becomes $3\pi/4$ and that of $w_-$
becomes $-5\pi/4$, which is outside the standard interval $(-\pi,\pi]$
of the phase of $w$ in \pref{wp:intw}; that is, $w$ is outside the 
standard Riemann sheet. In Figure 3 the path for the case $t=-0.99$ is  
shown.
The technical details will be worked out when writing the 
numerical algorithms.

\section{Concluding remarks}
\label{sec:conc}

In a future paper we will discuss the numerical aspects and describe
computer algorithms based on the integral representations given in this
paper.  Several quantities have to be calculated with great care.  For
example, straightforward use of $\psi(\theta)$ defined in \pref{pp:psi}
when $\theta$ is small, that is, at the saddle point, will give
cancellation of leading digits.  Also, to represent the functions for a
large range of the parameters scaling is needed.

When implementing the representations we will decide if the steepest
descent paths will be used or approximations of these paths, as we
suggested for the $W-$function in \S\,\ref{subsec:wwn}.  For example,
integrating in \pref{pp:intw} along the vertical line through the saddle
point $w_0$ gives a simpler representation than \pref{pp:inttheta}.
However, the integral along the vertical line has a non-real phase
function.  Another approximation of a steepest descent contour is given
in \pref{nm:vofu}.  We will investigate efficiency aspects in
combination with programming aspects in deciding which representation in
these examples should be used.

This also holds for the quite complicated steepest descent paths 
in \S\,\ref{subsubsec:wwpm}. We have not indicated in \pref{wpm:intv} the 
relation between $v$ and $u$ on the different parts of the path. This will be 
done during the implementation of the algorithms.

\section*{Acknowledgments}
\label{sec:ack}
 The authors thank
the referees for their careful reading of the manuscript 
and their valuable comments.
A. Gil acknowledges financial support from   Ministerio de 
 Ciencia y Tecnolog\'{\i}a                        
 (BFM2001-3878-C02-01).

\end{document}